\documentclass[12pt]{article}

 \newcommand{\beq}{\begin{equation}}
\newcommand{\eeq}{\end{equation}}
 \newcommand{\beth}{\begin{theo}}
\newcommand{\eth}{\end{theo}}

\makeatletter
\@addtoreset{equation}{section}
\makeatother

\newcommand{\Subset}{\subset\subset}
\newcommand{\supp}{{\rm supp}\,}
\newcommand{\Psux}{{ \Psi_{u,x}}}
\newcommand{\Psuo}{{ \Psi_{u,0}}}

\newcommand{\Psphy}{{ \Psi_{\varphi,y}}}
\newcommand{\psux}{{ \psi_{u,x}}}
\newcommand{\psuxp}{{ \psi_{u,x}^+}}

\newcommand{\okn}{1\le k\le n}
\newcommand{\LL}{{\cal L}}
\newcommand{\TT}{{\bf T}}
\newcommand{\DD}{{\bf D}}
\def\brv{{B_r(\varphi)}}
\def\srv{{S_r(\varphi)}}
\def\vp{\varphi}
\def\mrv{\mu_r^\varphi}

\def\mrV{\mu_r^\Phi}
\def\srV{{S_r(\Phi)}}
\newcommand{\aoP}{\gamma_1^\Phi}

\newcommand{\dpl}{{\cal D}_p^+}
\newcommand{\Ldpl}{{\cal LD}_p^+}
\newcommand{\Ldn}{{\cal LD}_{n-1}^+}
\newcommand{\Znp}{{ \bf Z}_+^n}
\newcommand{\Rn}{{ \bf R}^n}
\newcommand{\Rnm}{{ \bf R}_-^n}
\newcommand{\Rnp}{{ \bf R}_+^n}
\newcommand{\Cn}{{{ \bf C}^n}}
\newcommand{\Cp}{{ \bf C}^p}
\newcommand{\Cm}{{ \bf C}^m}
\newcommand{\Vol}{{\rm Vol}\,}

\newtheorem{theo}{Theorem}[section]
\newtheorem{cor}{Corollary}[section]

\newtheorem{prop}{Proposition}[section]

\begin{document}

\begin{center}
 {\bf \Large  Total Masses of Mixed Monge-Amp\`ere Currents}
\end{center}

\vskip0.5cm
\begin{center}
{\bf\large\sc Alexander Rashkovskii  }
\end{center}

\vskip1cm

  {\small {\sc Abstract.}  Mixed
Monge-Amp\`ere currents generated by plurisubharmonic functions
of logarithmic growth are studied. Upper bounds for their total
masses are obtained in terms of growth characteristics of the
functions. In particular, this gives a plurisubharmonic version of
D.~Bernstein's theorem on the number of zeros of polynomial
mappings in terms of the Newton polyhedra, and a representation
for Demailly's generalized degrees
of plurisubharmonic functions. }

\vskip 0.5cm

\section{Introduction}

Our starting point is the classical problem on numeric
characteristics for zero sets of polynomial mappings $P:\Cn \to
\Cm$. If $m\ge n$ and $P$ has discrete zeros, this is about the
total number of the zeros counted with multiplicities, and for
$m<n$ the characteristics are the projective volumes of the
corresponding holomorphic chains $Z_P$. When $m=1$, the volume
equals the degree of the polynomial $P$, while for $m>1$ the
situation becomes much more difficult. In particular, in the
general case no exact formulas can be obtained in terms of the
exponents and the problem reduces to finding appropriate upper
bounds. An example of such a bound is given by Bezout's theorem:
if $m= n$ and $P$ has discrete zeros, then their number does not
exceed the product of the degrees of the components of $P$. An
alternative estimate is due to Kouchnirenko \cite{Kou1},
\cite{Kou2}: the number of the zeros is at most $n!$ times the
volume of the {\it Newton polyhedron of $P$ at infinity} (the
convex hull of all the exponents of $P$ and the origin). A
refined version of the latter result was obtained by D.~Bernstein
\cite{Ber}. He showed that the number of (discrete) zeros of a
Laurent polynomial mapping $P$ on $({\bf C}\setminus\{0\})^n$ is
not greater than $n!$ times the mixed volume of the {\it Newton
polyhedra} (the convex hulls of the exponents) of the components
of $P$.

Here we put this problem into a wider context of pluripotential
theory. This can be done by considering plurisubharmonic functions
$u=\log|P|$ and studying the Monge-Amp\`ere operators
$(dd^cu)^p$; we use the notation $d=\partial + \bar\partial$,
$d^c= ( \partial -\bar\partial)/2\pi i$. The key relation is the
King-Demailly formula which implies that if the codimension of
the zero set is at least $p$ then $(dd^cu)^p\ge Z_P$ (with an
equality if $p=m\le n$). And the problem of estimating total
masses of the Monge-Amp\`ere operators of plurisubharmonic
functions $u$ of logarithmic growth was studied in \cite{R1}. In
particular, a relation was obtained in terms of the volume of a
certain convex set generated by the function $u$, which in case
$u=\log|P|$ is just the Newton polyhedron of $P$ at infinity.

On the other hand, the holomorphic chain $Z_P$ with $m=p\le n$ can
be represented as the wedge product of the currents (divisors)
$dd^c\log|P_k|$, $1\le k\le m$, which leads to consideration of
the {\it mixed} Monge-Amp\`ere operators
$dd^cu_1\wedge\ldots\wedge dd^cu_m$ and estimating their total
masses. Another motivation for this problem are generalized
degrees $\int_\Cn T\wedge (dd^c\vp)^p$ of positive closed currents
$T$ with respect to plurisubharmonic weights $\vp$  due to
Demailly \cite{D10}.

So, our main subject is mixed Monge-Amp\`ere currents generated by
arbitrary plurisubharmonic functions of logarithmic growth. Using
the approach developed in \cite{R1} we  get effective bounds for
the masses of the currents. As a consequence, this gives us a
plurisubharmonic version of Bernstein's theorem adapted, in
particular, for polynomial mappings of $\Cn$. In addition, we get
a representation for the generalized degrees of $(1,1)$-currents.

The paper can be considered as a continuation of \cite{R1}, and we
provide the reader with necessary definitions and statements from
there.

\section{Preliminaries and description of
results}\label{sect:prelim}

 We consider plurisubharmonic functions $u$ of logarithmic
growth in $\Cn$,
$$u(z)\le C_1\log^+|z|+C_2 $$
with some constants
$C_j=C_j(u)$. The collection of all such functions will be
denoted by $\LL(\Cn)$ or simply by $\LL$. Various results on such
functions are presented, e.\,g., in \cite{BT1}, \cite{BT2},
\cite{Le4}, \cite{Le5}, \cite{Si1}. For general properties of
plurisubharmonic functions and the complex Monge-Amp\`ere
operators, we refer the reader to \cite{Le2}, \cite{LeG},
\cite{Kl}, and \cite{D1}.

The {\it (logarithmic) type} of a function $u\in\LL$ is defined as
$$
\sigma(u)=\limsup_{|z|\to\infty}\frac{u(z)}{\log|z|}
$$
which can be viewed as the Lelong number of $u$ at infinity. The
corresponding counterpart for the directional (refined) Lelong
numbers are directional types \beq
\sigma(u,a)=\limsup_{z\to\infty} {u(z)\over \varphi_a(z)}
\label{eq:dirtype} \eeq with \beq
\varphi_a(z)=\sup_k\,a_k^{-1}\log|z_k|, \quad
a=(a_1,\ldots,a_n)\in\Rnp, \label{eq:Sa} \eeq see \cite{R1}. One
can also consider the types $\sigma(u,\varphi)$ with respect to
arbitrary plurisubharmonic exhaustive functions $\varphi\in\LL$,
\beq \label{eq:gentype} \sigma(u,\varphi)=\limsup_{z\to\infty}
{u(z)\over \varphi_(z)}. \eeq One more characteristic is the {\it
logarithmic multitype} $(\sigma_1(u),\ldots,\sigma_n(u))$
\cite{Le5}: \beq \sigma_1(u)=\sup\,\{\tilde\sigma_1(u;z'):
z'\in{\bf C}^{n-1}\} \label{eq:mtype} \eeq where
$\tilde\sigma_1(u;z')$ is the logarithmic type of the function
$u_{1,z'}(z_1)=u(z_1,z')\in\LL({\bf C})$ with $z'\in{\bf
C}^{n-1}$ fixed, and similarly for
$\sigma_2(u),\ldots,\sigma_n(u)$. For example, if $P$ is a
polynomial of degree $d_k$ in $z_k$, then $\sigma_k(\log|P|)=d_k$.

Another (and original) definition for the Lelong numbers is in
terms of the currents $dd^cu$, which works for arbitrary positive
closed currents. This leads to the notion of degree of a current.
 Let $\dpl(\Omega)$ be the collection of all closed
positive currents of bidimension $(p,p)$ on a domain
$\Omega\subset\Cn$. We will consider currents $T\in\dpl(\Cn)$
with finite projective mass, or {\it degree}
$$\delta(T)=\int_{\Cn}T\wedge
({1\over 2}dd^c\log(1+|z|^2))^p;$$ the set of all such currents is
denoted by $\Ldpl$. The  degree of $T\in\Ldpl$ can also be
represented as
$$\delta(T)=\int_{\Cn}T\wedge
(dd^c\log|z|)^p$$  and  as the density of the trace measure
$\sigma_T=T\wedge\beta_p$ of the current $T$:
$$
\delta(T)=\lim_{r\to\infty}\frac{\sigma_T(|z|<r)}{mes_{2p}(|z|<r)}.
$$
When $T=[A]$ is the current of integration over an algebraic set
$A$ of pure dimension $p$, the degree $\delta([A])$ coincides
with the degree of the set $A$ defined as the number of sheets in
the ramified covering map $A\to L$ to a generic $p$-codimensional
plane $L$. Note also that any current $T\in\Ldn$ has the form
$T=dd^cu$ with $u\in\LL$, and $\delta(dd^cu)=\sigma(u)$, see
\cite{LeG}.

The {\it generalized degrees} \beq \label{eq:gendeg}
\delta(T,\varphi):=\int_{\Cn}T\wedge (dd^c\varphi)^p \eeq with
respect to plurisubharmonic weights $\varphi$ were introduced in
\cite{D10} as a powerful tool for studying polynomial mappings
and algebraic sets.

 We are concerned with the problem of evaluation of
$$\mu(T,u_1,\ldots,u_p):=\int_{\Cn}T\wedge dd^cu_1\wedge\ldots\wedge dd^cu_p$$
for currents $T\in\Ldpl$ and functions $u_j\in\LL$ in terms of
the distribution of  $T$ and growth characteristics of  $u_k$.
The idea is to replace the functions $u_k$ by certain
plurisubharmonic functions $v_k$ with simpler asymptotic
properties. A relation between the corresponding total masses is
provided by Comparison Theorem~\ref{theo:comparison} which shows
that the value of $\mu(T,u_1,\ldots,u_p)$ is a function of the
asymptotic behavior of $u_k$ at infinity. The theorem is an
extension of B.A.~Taylor's theorem \cite{T} on the total mass of
$(dd^cu)^n$ of $u\in\LL\cap L_{loc}^\infty$. At the same time, it
is an analogue for Demailly's Second Comparison Theorem~5.9
\cite{D1} on generalized Lelong numbers.

By taking $v_k=\log|z|$ we get a bound in terms of the types of
$u_k$ and the degree of $T$ (Corollary~\ref{cor:finite})
$$
\mu(T,u_1,\ldots,u_p)\le \delta(T)\,\sigma(u_1)\ldots\sigma(u_p),
$$
and the choice $v_k=\vp_a$ leads to that in terms of the
corresponding directional characteristics
(Corollary~\ref{cor:dir}).

Sharper bounds are obtained with $v_k=\Psi_{u_k,x}$, the
indicators of $u_k$, introduced in \cite{R1} (see the definition
and basic properties in Section~\ref{sect:ind}). We have
$$
\mu(T,u_1,\ldots,u_p)\le\mu(T,\Psi_{u_1,x}^+,\ldots,\Psi_{u_p,x}^+)
$$
(Proposition \ref{prop:ind1}), and the problem is to find an
effective evaluation of the right-hand side. This can be done in
the case of $T=1$, $p=n$. Namely, the convex function
$\psuxp(t):=\Psi_{u_1,x}^+(\exp t_1,\ldots,\exp t_n)$, $t\in\Rn$,
is the support function to the convex set
$$
\Theta^u=\{a\in\Rn:\: \langle a,t\rangle\le\psuxp(t)\quad \forall
t\in\Rn\},
$$
and \beq\label{eq:mvol} \mu(\Psi_{u_1,x}^+,\ldots,\Psi_{u_n,x}^+)=
n!\,\Vol(\Theta^{u_1},\ldots,\Theta^{u_n}), \eeq Minkowski's mixed
volume of the sets $\Theta^{u_k}$ (Theorem~\ref{theo:indmain}).

The above considerations are applied in Section~\ref{sect:deg} to
investigation of the generalized degrees $\delta(T,\vp)$
(\ref{eq:gendeg}). By Proposition~\ref{prop:ind1}, we are reduced
to the values $\delta(T,\Psphy)$. In the particular situation of
$T=dd^cu$  we prove the relation
$\delta(dd^cu,\Psphy)=\delta(dd^c\Psux,\Psphy)$ for all
$x,y\in\Cn$. We study the so-called swept out Monge-Amp\`ere
measures of indicator weights in Theorem~\ref{theo:swept}. As a
consequence, we get a representation for $\delta(dd^cu,\Psphy)$ in
terms of the sets $\Theta^u$ and $\Theta^\vp$ in
Corollary~\ref{cor:swept}.

Finally, in Section~\ref{sect:alg} we specify our results for
currents generated by polynomial mappings. In particular, we
observe that (\ref{eq:mvol}) implies the following analogue for
D.~Bernstein's inequality (Corollary~\ref{cor:newtonmain}): the
projective volume
 $\delta(Z_P)$ of the holomorphic chain $Z_P$ generated by a
polynomial mapping $P=(P_1,\ldots,P_p)$ in general position,
$1\le p\le n$, has the bound
 $$\delta(Z_P)\le
 n!\,\Vol(G_1^+,\ldots,G_p^+,\Delta,\ldots,\Delta)$$
with $G_j^+$ the Newton polyhedron of the polynomial $P_j$ at
infinity and $\Delta=\{t\in\Rnp: \sum t_j\le 1\}$ the standard
simplex in $\Rn$. We also derive a number of other bounds (like
Bezout's and Tsikh's theorems) as direct consequences of our
general results on mixed Monge-Amp\`ere operators.

\section{Comparison theorem for mixed
operators}\label{sect:comparison}

A $q$-tuple of plurisubharmonic functions $u_1,\ldots,u_q$ will
be said to be {\it properly intersected}, or {\it in general
position}, with respect to a current $T\in\dpl$, $p\ge q$, if
their unboundedness loci $A_1,\ldots,A_p$ satisfy the condition:
for all choices of indices $j_1<\ldots<j_k,\ k\le q$, the
$(2q-2k+1)$-dimensional Hausdorff measure of the set $A_{j_1}\cap
\ldots\cap A_{j_k}\cap\supp T$ equals zero. If this is a case,
the current $T\wedge dd^cu_1\wedge\ldots\wedge dd^cu_q$ is
well-defined  and has locally finite mass (\cite{D1}, Th.~2.5).

We recall that a function $u$ in $\Cn$ is called {\it
semi-exhaustive} on a set $A$ if $\{u<R\}\cap A\Subset\Cn$ for
some real $R$, and {\it exhaustive} if this is valid for all $R$.

\beth \label{theo:comparison} {\rm (Comparison Theorem)} Let
$T\in\Ldpl$, $u_1,\ldots,u_p\in\LL$ be properly intersected with
respect to $T$, and $v_1,\ldots,v_p\in\LL$ be semi-exhaustive on
$\supp T$. If for any $\eta>0$,
$$
\limsup_{|z|\to\infty,\, z\in\supp T}
\frac{u_j(z)}{v_j(z)+\eta\log|z|}\le l_j,\quad 1\le j\le p,
$$
then $\mu(T,u_1,\ldots,u_p)\le l_1\ldots
l_p\,\mu(T,v_1,\ldots,v_p)$. \eth

{\it Proof.} It suffices to show that the conditions \beq
\label{eq:lim1} \limsup_{|z|\to\infty,\, z\in\supp T}
\frac{u_j(z)}{v_j(z)+\eta\log|z|} < 1,\quad 1\le j\le p, \
\forall\eta>0,\eeq imply \beq\label{eq:comp1}
\mu(T,u_1,\ldots,u_p)\le \mu(T,v_1,\ldots,v_p).\eeq Without  loss
of generality we can also take $\vp>0$ on $\Cn$.

For any $N>0$, the functions $u_{j,N}:=\max\{u_j,-N\}$ still
satisfy (\ref{eq:lim1}). Then for any $\eta>0$ and $C>0$, the set
$$
E_j(C)=\{z\in\supp T:\: v_j(z)+\eta\log|z|-C<u_{j,N}(z)\}
$$
is compactly supported in the ball $B_{\alpha_j}$ for some
$\alpha_j=\alpha_j(C,\eta,u_{j,N},v_j)$. Put
$\alpha=\max_j\alpha_j$, $E(C)=\cup_j E_j(C)$,
$F(C)=\cap_jE_j(C)$, $$w_{j,C}=\max\{v_j(z)+\eta\log|z|-C,
u_{j,N}\}.$$ Since $w_{j,C}=v_j(z)+\eta\log|z|-C$ near $\partial
B_\alpha\cap\supp T$, we have
\begin{eqnarray*}
\int_{B_\alpha}T\bigwedge_{1\le j\le p}dd^cw_{j,C}
&=&\int_{B_\alpha}T\bigwedge_{1\le j\le p}dd^c(v_j+\eta\log|z|)\\
&\le &\int_{\Cn}T\bigwedge_{1\le j\le p}dd^c(v_j+\eta\log|z|).
\end{eqnarray*}
Note that for any compact subset $K$ of $\supp T$ one can find
$C_K$ such that $K\subset F(C)$ for all $C>C_K$, so
$$
\int_{B_R}T\bigwedge_{1\le j\le p}dd^cw_{j,C} \le
\int_{\Cn}T\bigwedge_{1\le j\le p}dd^c(v_j+\eta\log|z|) $$
 for any
$R>0$ and all $C>C_R$. In addition,
$$
T\bigwedge_{1\le j\le p}dd^cw_{j,C}\to T\bigwedge_{1\le j\le
p}dd^cu_{j,N}
$$
as $C\to +\infty$ (the functions $w_{j,C}$ decrease to $u_{j,N}$)
and hence
\begin{eqnarray*}
\int_{B_R}T\bigwedge_{1\le j\le p}dd^cu_{j,N} &\le &
\limsup_{C\to\infty} \int_{B_R}T\bigwedge_{1\le j\le
p}dd^cw_{j,C}\\ &\le &  \int_{\Cn}T\bigwedge_{1\le j\le
p}dd^c(v_j+\eta\log|z|).
\end{eqnarray*}

Since $\eta$ is arbitrary, we derive the inequality
$$
\int_{B_R}T\bigwedge_{1\le j\le p}dd^cu_{j,N}\le
\int_{\Cn}T\bigwedge_{1\le j\le p}dd^cv_j
$$
and finally, letting $N\to\infty$,
$$
\int_{B_R}T\bigwedge_{1\le j\le p}dd^cu_{j}\le
\int_{\Cn}T\bigwedge_{1\le j\le p}dd^cv_j
$$
for any $R>0$, which gives us (\ref{eq:comp1}) and thus completes
the proof.

\medskip
{\it Remark.} As follows from the theorem,
$$
\mu(T,u_1,\ldots,u_p)\le\mu(T,\max\{u_1,\alpha_1\},\ldots,\max\{u_p,\alpha_p\})
$$
for any $\alpha\in{ \bf R}^p$, and the right-hand side is
independent of $\alpha$.
 The inequality here can be strict, which
follows from the consideration of the function
$u(z_1,z_2)=\log(|z_1|^2+|z_1z_2+1|^2)$. We have $\mu(u,u)=0$
while $\mu(u^+,u^+)=4$, the latter relation being verified by
comparing $u^+$ with the function
$\max\{\log|z_1|,\log|z_1z_2|,0\}$ whose total Monge-Amp\`ere
mass can be calculated by Proposition~\ref{prop:ind3} below. This
shows that the condition on the functions $v_j$ in
Theorem~\ref{theo:comparison} to be semi-exhaustive is essential.

\medskip
An immediate application of Comparison Theorem is the following
bound for the total mass of the current $T\wedge
dd^cu_1\wedge\ldots\wedge dd^cu_p$ in terms of the degree
$\delta(T)$ of $T$ and logarithmic types $\sigma(u_j)$ of $u_j$.

\begin{cor}\label{cor:finite}
If $T\in\Ldpl$ and $u_1,\ldots,u_p\in\LL$ are properly intersected
with respect to $T$, then
$$
\mu(T,u_1,\ldots,u_p)\le \delta(T)\,\sigma(u_1)\ldots\sigma(u_p).
$$
In particular, $T\wedge dd^cu_1\wedge\ldots\wedge dd^cu_{p-k}\in
{\cal LD}_k^+$, $0\le k\le p$.
\end{cor}

Moreover, we have a refined bound in terms of the generalized
characteristics $\sigma(u,\vp)$ and $\delta(T,\varphi)$ with
respect to weights $\vp$ (\ref{eq:gentype}), (\ref{eq:gendeg}).

\begin{cor}\label{cor:dir}
Let $T,u_1,\ldots,u_p$ satisfy the conditions of
Corollary~\ref{cor:finite} and $\vp\in\LL$ be an exhaustive
weight. Then
$$
\mu(T,u_1,\ldots,u_p)\le
\delta(T,\varphi)\,\sigma(u_1,\vp)\ldots\sigma(u_p,\vp).
$$
\end{cor}

Being specified to the case $T=1$, $p=n$ and $\vp=\vp_a$ this
gives us the following bound in terms of the directional types
$\sigma(u_j,a)$ of $u_j$ (\ref{eq:dirtype}).

\begin{cor}
If the functions $u_1,\ldots,u_n\in\LL$ are properly intersected,
then
$$
\mu(u_1,\ldots,u_n)\le
\inf_{a\in\Rnp}\frac{\sigma(u_1,a)\ldots\sigma(u_n,a)}{a_1\ldots
a_n}.
$$
\end{cor}

\section{Bounds in terms of indicators}\label{sect:ind}

More precise bounds can be obtained by means of indicators of
functions from the class $\LL$.

Developing the notion of local indicator introduced in \cite{LeR},
the (global) {\it indicator} of a function $u\in\LL$  at
$x\in\Cn$ was defined in \cite{R1} as
$$\Psux(y)=\lim_{R\to +\infty}R^{-1}\sup\{u(z):\: |z_k-x_k|\le
 |y_k|^R,\ \okn\} \ {\rm for}\ y_1\ldots y_n\neq 0,$$
and it extends to  a plurisubharmonic function of the class
$\LL$. The indicator controls the behavior of $u$ in the whole
$\Cn$ (\cite{R1}, Th.~1): \beq\label{eq:indgrowth}
 u(z)\le\Psux(z-x)+C_x\quad\forall z\in\Cn
\eeq with $C_x$ the supremum of $u$ on the unit polydisk centered
at $x$. Besides, it is a (unique) logarithmic tangent to $u$ at
$x$, i.\,e., the weak limit in $L_{loc}^1(\Cn)$ of the functions
\beq\label{eq:u_m}u_m(y)= m^{-1} u(x_1+y_1^m,\ldots,
x_n+y_n^m)\eeq as $m\to\infty$ (\cite{R1}, Th.~2).

Note that the indicator of $\log|z|$ at $ x$ equals $\max_k
\log|y_k|$ if $x=0$, and $\max_k \log^+|y_k|$ for any other point
$x$.

The asymptotic characteristics (types) of $u$ can be easily
expressed in terms of its indicator or, which is more convenient,
of the convex image
$$\psux(t)=\Psux(e^{t_1},\ldots,e^{t_n}),\quad t\in\Rn,$$
of the
indicator (\cite{R1}, Prop.~3): for each $x\in\Cn$, 
$$\sigma(u)=\sigma(u,(1,\ldots,1))=\psux(1,\ldots,1), $$ 
\beq\label{eq:sigmaa} \sigma(u,a)=\psux(a)\quad\forall a\in\Rnp,
 \eeq \beq\label{eq:sigmak}
\sigma_k(u)=\psux(e_k)\quad \okn\eeq with $e_1,\ldots,e_n$ the
standard basis of $\Rn$. Note that the restriction of $\psux$ to
$\Rnp$ is independent of $x\in\Cn$ (\cite{R1}, Prop.~7).

By (\ref{eq:indgrowth}), Theorem~\ref{theo:comparison} implies

\begin{prop}\label{prop:ind1}
Let $T\in\Ldpl$, $u_1,\ldots,u_p\in\LL$ be properly intersected
with respect to $T$, and $x\in\Cn$. Then
$$
\mu(T,u_1,\ldots,u_p)\le
\mu(T,\Psi_{u_1,x}^+,\ldots,\Psi_{u_p,x}^+).
$$
\end{prop}

\begin{cor}\label{cor:ind2}
Let $u_1,\ldots,u_n\in\LL$ be properly intersected and
$x^k\in\Cn$, $\okn$. Then
\begin{eqnarray*}
\mu(u_1,\ldots,u_n) &\le &
\mu(u_1,\ldots,u_{n-1},\Psi_{u_n,x^n}^+)\le
\mu(u_1,\ldots,\Psi_{u_{n-1},x^{n-1}}^+,\Psi_{u_n,x^n}^+)\ldots\nonumber\\
&\le & \mu(u_1,\Psi_{u_1,x^1}^+,\ldots,\Psi_{u_n,x^n}^+) \le
\mu(\Psi_{u_1,x^1}^+,\ldots,\Psi_{u_n,x^n}^+).
\end{eqnarray*}
\end{cor}

{\it Remark.} The choice of $x\in\Cn$ can affect the value of the
total Monge-Amp\`ere mass of the indicators. For example, let
$u(z_1,z_2)={1\over 2}\log(1+|z_1z_2|^2)$, then
$\Psi_{u,0}(y)=\log^+|y_1y_2|$ has zero mass, while the mass of
$\Psi_{u,(1,1)}(y)=\max\{\log^+|y_1|,\log^+|y_2|,\log^+|y_1y_2|\}$
equals $2$ (see Proposition~\ref{prop:ind3} below).

\medskip

To get an interpretation for the masses of indicators, we proceed
as in \cite{R1}. Let $\Phi$ be an abstract indicator in $\Cn$,
that is, a plurisubharmonic function  depending only on
$|z_1|,\ldots,|z_n|$ and satisfying the homogeneity condition
\beq\label{eq:homo}
\Phi(|z_1|^c,\ldots,|z_n|^c)=c\,\Phi(|z_1|,\ldots,|z_n|)\quad\forall
c>0. \eeq The functions $\vp_a$ (\ref{eq:Sa}) are particular
examples of the indicators.
 It is clear that $\Phi\le 0$ in the unit polydisk $\DD$ and
is strictly positive on
$$
\DD^{-1}=\{z\in\Cn:\; |z_k|>1,\ \okn\}
$$
unless $\Phi\equiv 0$.

Let us assume $\Phi\ge 0$ on $\Cn$. In this case, $(dd^c\Phi)^n$
is supported by the distinguished boundary $\TT$ of $\DD$
(\cite{R1}, Th. 6).
 Denote
\beq\label{eq:vp}
 \vp(t)=\Phi(e^{t_1},\ldots,e^{t_n}),\quad
t\in\Rn, \eeq the convex image of $\Phi$ in $\Rn$, and
\beq\label{eq:Theta} \Theta^\Phi=\{a\in\Rn:\: \langle
a,t\rangle\le\vp(t)\quad \forall t\in\Rn\}. \eeq
 It is easy to see that $\Theta^\Phi$ is a convex
compact subset of $\overline{\Rnp}$. By the construction, $\vp$ is
the support function of $\Theta^\Phi$. The {\it real}
Monge-Amp\'ere operator applied to $\vp$ gives us the
$\delta$-function $\delta_0$ with the mass $\Vol(\Theta^\Phi)$.
By comparing the real and complex Monge-Amp\'ere operators we get
\beq\label{eq:muphi} \mu(\Phi,\ldots,\Phi)=n!\,\Vol(\Theta^\Phi)
\eeq (see the details in \cite{R1}, Th.~6).

This can be extended to the mixed Monge-Amp\'ere operators of
indicators:

\begin{prop}\label{prop:ind3}
 Let
$\Phi_1,\ldots,\Phi_n$ be nonnegative indicators. Then
$$dd^c\Phi_1\wedge\ldots\wedge dd^c\Phi_n=\mu\,dm$$ with $dm$ the
normalized Lebesgue measure on $\TT$, and
$$
\mu=\mu(\Phi_1,\ldots,\Phi_n)=n!\,\Vol(\Theta^{\Phi_1},\ldots,\Theta^{\Phi_n}).
$$
\end{prop}

{\it Proof.} By the polarization formula for the complex
Monge-Amp\'ere operator, \beq\label{eq:polar}
\bigwedge_kdd^c\Phi_k=\frac{(-1)^n}{n!}\sum_{j=1}^n(-1)^j
\sum_{1\le i_1<\ldots<i_j\le n}\left(dd^c\sum_{k=1}^j \Phi_{j_k}
\right)^n. \eeq

Since the sum of indicators is an indicator as well, the support
of $ \wedge_k dd^c\Phi_k$ is a subset of $\TT$. In view of the
translation invariance of this measure, it has the form $\mu\,
(2\pi)^{-n}\,d\theta_1\ldots d\theta_n$.

By (\ref{eq:muphi}), each term of the right-hand side of
(\ref{eq:polar}) equals $n!$ times the corresponding volume, so
the left-hand side is just the mixed volume of the convex sets
$\Theta^{\Phi_k}$.

\medskip
Since $\Phi=\Psi_{u,x}^+$ satisfies (\ref{eq:homo}),
Corollary~\ref{cor:ind2} and Proposition \ref{prop:ind3} give us
the main result of the section.

\beth\label{theo:indmain} Let functions $u_1,\ldots,u_n\in\LL$ be
properly intersected, $x^k\in\Cn$, $\okn$. Then
$$
\mu(u_1,\ldots,u_n) \le
n!\,\Vol(\Theta^{\Phi_1},\ldots,\Theta^{\Phi_n})
$$
with the sets $\Theta^{\Phi_k}$ defined by (\ref{eq:vp}) -
(\ref{eq:Theta}) for  $\Phi=\Phi_k=\Psi_{u_k,x^k}^+$. \eth

As a consequence we can get a bound for $\mu(u_1,\ldots,u_n)$ in
terms of the types $\sigma_k(u_j)$ of $u_j$ with respect to $z_k$
(\ref{eq:mtype}).

\begin{cor}\label{cor:perm}
If $u_1,\ldots,u_n\in\LL$ are properly intersected, then
$$\mu(u_1,\ldots,u_n) \le n!\,{\rm
per}\left(\sigma_k(u_j)\right)_{j,k=1}^n,$$ where ${\rm per}\,A$
denotes the permanent of the matrix $A$.
\end{cor}

{\it Proof.} As follows from (\ref{eq:sigmak}), the set
$\Theta^{\Phi_j}$ is a subset of the rectangle
$[0,\sigma_1(u_j)]\times\ldots \times [0,\sigma_n(u_j)]$, and the
mixed volume of the rectangles $[0,a_{1j}]\times\ldots \times
[0,a_{nj}]$, $1\le j\le n$, equals ${\rm
per}\left(a_{jk}\right)_{j,k=1}^n$.

\section{Degrees with respect to plurisubharmonic
weights}\label{sect:deg}

Given a subset $A$ of $\Cn$, we denote by $W(A)$ the collection
of all continuous, as mappings to $[-\infty,+\infty)$, functions
({\it weights}) $\varphi\in\LL$ exhaustive on $A$.

Let $T\in\Ldpl$,  then the current (measure)
$T\wedge(dd^c\varphi)^p\in{\cal LD}_0^+$ is well-defined for any
$\vp\in W(\supp T)$. We put
$$
\delta(T,\varphi)=\int_{\Cn}T\wedge(dd^c\varphi)^p
$$
be the (generalized) degree of $T$ with respect to the weight
$\varphi$, see \cite{D10}.

Observe that $\delta(T,\vp)=\delta(T,\vp^+)$ since $\vp$ is
assumed to be exhaustive on $\supp T$. Note also that
$\delta(T,\varphi)=\delta(T)$ if $\varphi(z)=\log|z|$.

Generalized degrees of currents can be viewed as generalized
Lelong numbers at infinity, and we start here with two
semicontinuity  properties parallel to those for the Lelong
numbers, cf. \cite{D1}.

\begin{prop}
\label{prop:sem1} Let $T_m,T\in \Ldpl$, $T_m\to T$. Then for any
weight $\varphi$ from $ W\left(\overline{\cup_m\supp T_m}\right)$,
$$
\delta(T,\varphi)\le\liminf_{m\to\infty}\delta(T_m,\varphi).
$$
\end{prop}

{\it Proof} follows immediately from Prop.~3.12 of \cite{D1}.

\begin{prop}\label{prop:sem2}
Let weights $\varphi_k,\,\varphi\in W(\supp T)$ be such that for
some $t\in{\bf R}$ the functions $\max\{\varphi_k,t\}$ converge to
$\max\{\varphi,t\}$ uniformly on compact subsets of $\Cn$. Then
$$
\delta(T,\varphi)\le\liminf_{m\to\infty}\delta(T,\varphi_k).
$$
\end{prop}

{\it Proof.} Since $\delta(T,\max\{\varphi,t\})=\delta(T,\vp)$
for any weight $\vp\in W(\supp T)$, we can take $\vp_k\to\vp$
uniformly on compact subsets of $\Cn$.

For any $R>0$ consider $\eta\in C^\infty(\Cn)$, $0\le \eta\le 1$,
such that $\supp\eta\subset B_R$ and $\eta\equiv 1$ on $B_{R/2}$.
The relation
$$
\lim_{k\to\infty}\int\eta T\wedge(dd^c\vp_k)^p= \int\eta
T\wedge(dd^c\vp)^p
$$
implies
$$
\liminf_{k\to\infty}\delta(T,\vp_k)\ge \int\eta
T\wedge(dd^c\vp)^p,
$$
and the assertion follows.

\medskip

Comparison Theorem~\ref{theo:comparison} for the degrees reads as
follows.

\begin{prop}\label{prop:compwei}
If two weights $\varphi,\psi\in W(\supp T)$ for a current
$T\in\Ldpl$ and
$$
\limsup_{|z|\to\infty,\, z\in\supp T} \frac{\vp(z)}{\psi(z)}\le l,
$$
then $\delta(T,\varphi)\le l^p\,\delta(T,\psi)$.
\end{prop}

Being applied to indicators, this gives us

\begin{cor}\label{cor:curindwei}
For any current $T\in\Ldpl$, any weight $\vp\in W(\supp T)$ and
$y\in\Cn$,
$\delta(T,\varphi)\le\delta(T,\Psi_{\vp,y}^+)=\delta(T,\Psphy)$.
\end{cor}

More can be said in the case of $T=dd^cu$, $u\in\LL$. A weight
$\Phi$ will be called {\it homogeneous} if it depends only on
$|z_1|,\ldots,|z_n|$ and  satisfies the homogeneity condition
(\ref{eq:homo}). In other words, homogeneous weights are
(abstract) indicators. Note that a homogeneous weight $\Phi$ is
exhaustive on $\Cn$ if and only if $\Phi>0$ on $\Cn\setminus \DD$.
In addition, it is easy to see that the type $\sigma(u,\Phi)$
with respect to such a weight can be computed as
$$
\sigma(u,\Phi)=\max\{\Psuo(z):\: \Phi(z)=1\}.
$$
It is interesting that the degrees $\delta(dd^cu,\Phi)$ can also
be represented in terms of the indicators.

\beth\label{theo:indwei} For any function $u\in\LL$, any $x\in\Cn$
and any homogeneous weight $\Phi\in W(\Cn)$, the equality
$\delta(dd^cu,\Phi)=\delta(dd^c\Psux,\Phi)$ holds. In particular,
$\delta(dd^cu,\varphi)\le\delta(dd^cu,\Psphy)=\delta(dd^c\Psux,\Psphy)$
for any  weight $\varphi\in W(\Cn)$ and $y\in\Cn$. \eth

{\it Proof.} Consider the family of functions $u_m$ defined by
(\ref{eq:u_m}), $x\in\Cn$.  As was mentioned in
Section~\ref{sect:ind}, $u_m$ converge (in $L_{loc}^1$) to $\Psux$
as $m\to\infty$, so $dd^cu_m\to dd^c\Psux$. By
Proposition~\ref{prop:sem1},
$$
\delta(dd^c\Psux,\Phi)\le\liminf_{m\to\infty}\delta(dd^c u_m,
\Phi).
$$
However the homogeneity of $\Phi$ gives us $\delta(dd^c u_m,
\Phi)=\delta(dd^c u, \Phi)$ for each $m$, so
$\delta(dd^cu,\Phi)\ge\delta(dd^c\Psux,\Phi)$ and the desired
equality follows then from Corollary~\ref{cor:curindwei}. The
theorem is proved.

\medskip

{\it Remark.} A relation between $\sigma(u,\Phi)$ and
$\delta(dd^cu,\Phi)$ will be given in Corollary~\ref{cor:swept}.

\medskip

The generalized degrees can be represented by means of the swept
out Monge-Amp\`ere measures introduced by Demailly \cite{D3}. For
$\vp\in W(\Cn)$, let $\brv=\{z:\vp(z)<r\}$,
$\srv=\{z:\vp(z)=r\}$, $\vp_r=\max\,\{\vp,r\}$. The {\it swept
out Monge-Amp\`ere measure} $\mrv$ is defined as
$$
\mrv
=\left(dd^c\vp_r\right)^n
-\chi_r (dd^c\vp)^n
$$
with $\chi_r$ the characteristic function of $\Cn\setminus
B_r(\vp)$. It is a positive measure on $\srv$ with the total mass
$\mrv(\srv)=(dd^c\vp)^n(\brv)$. If $\supp(dd^c\vp)^n\subset
B_R(\vp)$, then $\mrv=(dd^c\vp_r)^n$ for all $r>R$.

 For $\vp(z)=\log|z-x|$, $\mrv$ is the normalized Lebesgue measure
on the sphere $\{z:|z-x|=e^r\}$, and for $\vp=\vp_a$ it is
supported by the set \beq \label{eq:tra} \TT_{ra}=\{z:\: z_k=\exp
( ra_k+i\theta_k),\okn\}\eeq and has the form $\mu_r^{\vp_a}=
(a_1\ldots a_n)^{-1}(2\pi)^{-n}\,d\theta_1 \ldots d\theta_n$
\cite{D1}.

The role of the measures $\mrv$ is clarified by the
Lelong-Jensen-Demailly formula \cite{D3}, \cite{D1}: for any
function $u$ plurisubharmonic in $B_R(\vp)$,
$$
\mrv(u)-\int_\brv u(dd^c\vp)^n=\int_{-\infty}^r\int_{B_t(\vp)}
dd^cu\wedge(dd^c\varphi)^{n-1}\,dt\quad \forall r<R.
$$

\beth (cf. \cite{D3}, \cite{D1}) \label{prop:asdir} Let $u\in\LL$
and $\vp\in W(\Cn)$. Then \beq\label{eq:unb}
\limsup_{r\to\infty}{\mrv(u)\over r}\le\delta(dd^cu,\vp)\le
\liminf_{r\to\infty}{\mrv(u^+)\over r}. \eeq
 If, in addition,
 \beq\label{eq:bounded}
 \supp(dd^c\vp)^n\subset B_{r_0}(\vp)
 \eeq
 for some $r_0$, then $r\mapsto\mrv(u)$ is a convex function of $r\in
(r_0,\infty)$ and \beq\label{eq:bounded1} \delta(dd^c
u,\vp)=\lim_{r\to +\infty}{\mrv(u)\over r}. \eeq \eth

{\it Proof.} From the Lelong-Jensen-Demailly formula we have for
any $r>r_0$, \beq\label{eq:LJD1}
 \mrv(u)=\int_{B_{r_0}(\vp)}
u(dd^c\vp)^n+\int_{\brv\setminus B_{r_0}(\vp)}
u(dd^c\vp)^n+\int_{-\infty}^r \delta(dd^cu,\vp,t)\,dt \eeq
 with
$$
\delta(dd^cu,\vp,t)=\int_{B_t(\vp)}
dd^cu\wedge(dd^c\varphi)^{n-1}.
$$
If $\vp$ satisfies (\ref{eq:bounded}) then the right-hand side is
a convex function of $r$, and (\ref{eq:bounded1}) follows. When
(\ref{eq:bounded}) is not assumed, take any $\epsilon>0$ and
choose $r_0$ such that $(dd^c\vp)^n(\Cn\setminus
B_{r_0}(\vp))<\epsilon$ and
$u(z)\le(\sigma(u,\vp)+\epsilon)\vp(z)$ for all $z\in\Cn\setminus
B_{r_0}(\vp)$. Then we get
$$
\mrv(u)\le Const + (\sigma(u,\vp)+\epsilon)r\epsilon
+\int_{-\infty}^r \delta(dd^cu,\vp,t)\,dt
$$
which gives us the first inequality in (\ref{eq:unb}). To get the
second inequality, consider the functions
$u_N(z)=\max\{u(z),-N\}$, $N>0$. By Proposition~\ref{prop:sem1},
the number $N$ can be chosen such that $\delta(dd^cu_N,\vp)\ge
\delta(u,\vp)+\epsilon$. Application of (\ref{eq:LJD1}) to the
function $u_N$ gives us
$$
\mrv(u^+)\ge\mrv(u_N)\ge Const -N\epsilon+\int_{-\infty}^r
\delta(dd^cu,\vp,t)\,dt
$$
and thus finishes the proof.

\medskip

The structure of the swept out Monge-Amp\`ere measures for
homogeneous weights  is given by

\beth\label{theo:swept} Let $\Phi\in W(\Cn)$ be a homogeneous
weight. For any function $u$ plurisubharmonic in $B_R(\Phi)$,
$R>0$, the swept out Monge-Amp\`ere measure $\mrV$ on the set
$\srV$, $0<r<R$, is determined by the formula
$$
\mrV(u)=n!\,\int_{E^\Phi}\lambda(u,rt)\,d\aoP(t)
$$
where $\lambda(u,rt)$ is the mean value of $u$ over the
distinguished boundary of the polydisk $\{|z_k|<\exp\{rt_k\},\
\okn\}$ and the measure $\aoP$ on the set $E^\Phi$ of extreme
points of the convex set $\{t\in\Rn: \Phi(e^{t_1},\ldots,
e^{t_n})\le 1\}$ is given by the relation $\aoP(F)={\rm
Vol}\,\Theta_F^\Phi$ for compact subsets $F$ of $E^\Phi$, the set
$\Theta_F^\Phi$ being defined by relations (\ref{eq:gf}),
(\ref{eq:gth}), and (\ref{eq:LPhi}) below. \eth

\begin{cor}\label{cor:swept}
For any $\varphi\in W(\Cn)$, $u\in\LL$ and $x,y\in\Cn$,
$$\delta(dd^cu,\vp)\le\delta(dd^cu,\Psi_{\vp,y})=
n!\,\int_{E^\Phi}\psi_{u,x}(t)\,d\aoP(t)$$ with
$\Phi=\Psi_{\vp,y}$ and $\psux$ the convex image of the indicator
$\Psux$. In particular, for any homogeneous weight $\Phi$,
$\delta(dd^cu,\Phi)\le\sigma(u,\Phi)\,\mu(\Phi,\ldots,\Phi)$.
\end{cor}

{\it Remark.} A description for the swept out Monge-Amp\`ere
measures for (negative) local indicators was given in \cite{R4}
and it resulted in the fact that the generalized Lelong numbers of
a plurisubharmonic function can be expressed in terms of those
with respect to the weights $\vp_a$ (its directional Lelong
numbers). As follows from Corollary~\ref{cor:swept}, this is no
longer true for the generalized weights.

\medskip


{\it Proof of Theorem \ref{theo:swept}.}  Since $(dd^c\Phi)^n=0$
on $\{\Phi>0\}$, we have $\mrV=(dd^c\Phi_r)^n$ for each $r>0$. By
the rotation invariance,
$$
\mrV=(2\pi)^{-n}d\theta\otimes d\rho_r^\Phi
$$
with some measure $\rho_r^\Phi$ supported by $S_r(\Phi)\cap\Rn$.
Moreover, since $\mrV$ has no masses on the pluripolar set
$S_r(\Phi)\cap\{z:z_1\ldots z_n=0\}$, we can pass to the
coordinates $z_k=\exp\{t_k+i\theta_k\}$, $-\infty<t_k<\infty$,
$\okn$. The functions
$$
\vp(t):=\Phi(e^{t_1},\ldots,e^{t_n}),\quad
\vp_r(t)=\Phi_r(e^{t_1},\ldots,e^{t_n})= \max\,\{\vp(t),r\}
$$
are convex in $\Rn$ and increasing in each $t_k$. Simple
calculations show that in these coordinates $\rho_r^\Phi$
transforms into the measure
$$
\gamma_r^{\Phi}=n!\,{\cal MA}[\vp_r]
$$
where ${\cal MA}$ is the {\it real} Monge-Amp\`ere operator, see
the details, e.g., in \cite{R}. We recall that for smooth
functions $v$,
$$
{\cal MA}[v]=\det\left({\partial^2 v\over\partial t_j\partial t_k}
\right)\,dt,
$$
and it can be extended as a positive measure to any convex
function (see \cite{RaT}). So, we have
\begin{eqnarray*}
\mrV(u) &=& \int_{\Rn}(2\pi)^{-n}\int_{[0,2\pi]^n}
u(z_1e^{i\theta_1},
     \ldots, z_ne^{i\theta_n})\,d\theta\,d\rho_r^\Phi(|z_1|,\ldots,|z_n|)\\
&=&n!\,\int_{\Rn}\lambda(u,t)\,d\gamma_r^{\Phi}(t)
=n!\,\int_{\Rn}\lambda(u,rt)\,d\aoP(t)
\end{eqnarray*}
since $\vp_r(t)=r\vp_{1}(t/r)$,
and we only need to find an explicit expression for the measure
$\aoP$  supported in the level set \beq\label{eq:LPhi}
 L^\Phi=\{t\in\Rn:
\vp(t)=1\}. \eeq

 As follows from properties of the real Monge-Amp\`ere operator,
  \beq \int_{F} {\cal MA}[\vp_1]={\rm
Vol}\,(\omega(F,\vp_1))\quad\forall F\subset L^\Phi, \label{eq:u3}
\eeq where
$$
\omega(F,\vp_1)=\bigcup_{t^0\in F} \{a\in\Rn:\: \vp_1(t)\ge 1+
\langle a,t-t^0\rangle\quad\forall t\in \Rn\}
$$
is the gradient image of the set $F$.

Given a subset $F$ of $L^\Phi$, we put \beq
\Gamma_F^\Phi=\{a\in\Rnp:\: \sup_{t\in F}\langle a,t\rangle =
\sup_{t\in L^\Phi}\langle a,t\rangle=1\} \label{eq:gf} \eeq and
\beq \Theta_F^\Phi=\{\lambda a:\: 0\le\lambda\le 1,\ a\in
\Gamma_F^\Phi\}. \label{eq:gth} \eeq Note that
$\Theta_{L^\Phi}^\Phi$ is a bounded convex subset of $\Rnp$ and
$\vp$ is  its support function. We claim that {\sl for any compact
subset $F$ of} $L^\Phi$, $\Theta_F^\Phi=\omega(F,\vp_1)$.

 If $a\in\omega(F,\vp_1)$ then for some $t^0\in F$,
\beq \langle a,t^0\rangle\ge\langle
a,t\rangle-\vp_1(t)+1\quad\forall t\in\Rn. \label{eq:th1} \eeq In
particular, \beq \langle a,t^0\rangle\ge\langle
a,t\rangle\quad\forall t\in L^\Phi. \label{eq:th2} \eeq When
$t=c\,t_0$ with $c>1$, (\ref{eq:th1}) implies $\langle
a,t^0\rangle\le 1$. In view of (\ref{eq:th2}) it means that
$a\in\Theta_F^\Phi$ and thus
$\omega(F,\vp_1)\subset\Theta_F^\Phi$.

Let now $a\in\Theta_F^\Phi$, so $a=\lambda a^0,\
a^0\in\Gamma_F^\Phi,\ 0\le\lambda\le 1$. Then there is a point
$t^0\in F$ such that
$$
\langle a,t^0\rangle=\sup_{t\in F}\langle a,t\rangle = \sup_{t\in
L^\Phi}\langle a,t\rangle=\lambda.
$$
Take any $t\in\Rn$. If $\vp(t)\le 1$, then $\langle
a,t^0\rangle\ge\langle a,t\rangle$ and thus $\vp_1(t)\ge 1
+\langle a,t-t^0\rangle$. If $\vp(t)=\alpha>1$, then $t/\alpha\in
L^\Phi$ and
\begin{eqnarray*}
\langle a,t\rangle-\vp_{1}(t)+1&=&\alpha\langle
a,t/\alpha\rangle+1
-\alpha\le\alpha\sup_{s\in L^\Phi}\langle a,s\rangle +1-\alpha\\
&=&\alpha\langle a,t^0\rangle+1-\alpha =\alpha\lambda+1-
\alpha\le\lambda=\langle a,t^0\rangle,
\end{eqnarray*}
so $a\in\omega(F,\vp_1)$. The claim is proved.

Finally, let $E^\Phi$ be the set of extreme points of $L^\Phi$
(that is, those not situated inside intervals on $L^\Phi$).
 As $L^\Phi\subset\{t\in\Rn\setminus\Rnm:\: t_k\le b_k,\ \okn\}$ for
 some $b\in\Rnp$, we have
$$
\sup_{t\in L^\Phi}\langle a,t\rangle = \sup_{t\in E^\Phi}\langle
a,t\rangle\quad\forall a\in\Rnp,
$$
so that $\Theta_{L^\Phi}^\Phi=\Theta_{E^\Phi}^\Phi$. Hence
$\aoP(L^\Phi)=\aoP(E^\Phi)$ and then $\supp\aoP\subset E^\Phi$.
The proof is complete.

\section{Algebraic case: Newton polyhedra}\label{sect:alg}
Here we test our results for the case of currents generated by
polynomial mappings.

 When $u=\log|P|$ for a polynomial $P$, $\sigma(u)$ is the
degree of $P$, $\sigma_k(u)$ is its degree with respect to $z_k$,
$\psux(t)=\max \,\{\langle t,J\rangle :\: J\in\omega_x(P)\}$ where
$$
\omega_x(P)=\{J\in\Znp:\: {\partial^JP\over\partial z^J}(x)\neq
0\},
$$
see \cite{R1}. Note that the maximum is attained on the set of
extreme points $E_x(P)$ of the set $\omega_x(P)$. In particular,
$E_0(P)$ coincides with the set of extreme points of the convex
hull  of the exponentials of the polynomial $P$.
This means that the set $\Theta^\Phi$ with $\Phi=\Psuo$ equals
\beq\label{eq:NP}G^+(P)= conv\,(E_0(P)\cup\{0\}),\eeq the {\it
Newton polyhedron  of $P$ at infinity} as defined in \cite{Kou1}.

Let now $P=(P_1,\ldots,P_p)$ be a polynomial mapping,
$u_j=\log|P_j|$, then $Z_j=dd^cu_j$ is the divisor of $P_j$. Let
the zero sets $A_j$ of $P_j$ be properly intersected, that is,
${\rm codim}\,A_{j_1}\cap\ldots\cap A_{j_m}\ge m$ for all choices
of indices $j_1,\ldots,j_m$, $m\le p$.  Then the holomorphic chain
$Z$ of the mapping $P$ is the intersection of the divisors $Z_j$:
$Z=Z_1\wedge\ldots\wedge Z_p$ (\cite{D1}, Prop.~2.12).

In this setting,
Corollary~\ref{cor:finite} turns to the bound
$$
\int_{\Cn}T\wedge Z\le \delta(T)\,\delta_1\ldots\delta_p
$$
via the degrees $\delta_j$ of $P_j$. For $T=1$ this gives
Bezout's inequality for the projective volume of the chain $Z$.
The specification of Corollary~\ref{cor:dir} with $a\in\Znp$
gives a bound in terms of the degrees in terms of a
quasi-homogeneous mapping $P(z^a)$, a global counterpart for
Tsikh-Yuzhakov's theorem on multiplicity of holomorphic mappings
in terms of the weighted initial homogeneous polynomial terms
\cite{YTs}. And Corollary~\ref{cor:perm} implies a result of Tsikh
\cite{Ts}.

Theorem~\ref{theo:indmain} now takes the form

\begin{cor}\label{cor:newtonmain}
The degree (projective volume) $\delta(Z)$ of the holomorphic
chain $Z$ generated by a polynomial mapping $P=(P_1,\ldots,P_p)$,
$p\le n$, with the zero sets of its components $P_k$ properly
intersected,
 has the bound
$$
\delta(Z)\le n!\,\Vol(G_1^+,\ldots,G_p^+,\Delta,\ldots,\Delta),
$$
$G_j^+$ being the Newton polyhedron of the polynomial $P_j$ at
infinity (defined by (\ref{eq:NP})) and $\Delta=\{t\in\Rnp: \sum
t_j\le 1\}$ the standard simplex in $\Rn$. In particular, when
$p=n$, the number of zeros of $P$ counted with their
multiplicities does not exceed $n!\,\Vol(G_1^+,\ldots,G_n^+)$.
\end{cor}

When $P_j(0)=0$, the set $G_j^+$ is strictly greater than the
convex hull $E_0(P_j)$ of the set $\omega_0(P_j)$ appearing in
Bernstein's theorem. But in return we take care of all the zeros
while Bernstein's theorem estimates only those in $({\bf
C}\setminus\{0\})^n$. And actually no bound for the total number
is possible in terms of just the convex hulls of the exponents
(see, e.g., $f(z)\equiv z$).

An algebraic specification of Theorem \ref{theo:swept} and
Corollary~\ref{cor:swept} is as follows. Let
 $\Phi$ be the indicator of $\log|P|$ for a polynomial mapping
$P:\Cn\to\Cp$, $p\ge n$, with discrete zeros. Then the set
$\Gamma_{L^\Phi}^\Phi$ is the Newton diagram for $P$ and
$\Theta_{L^\Phi}^\Phi=G^+(P)$ is the Newton polyhedron for $P$ at
infinity, that is, the convex hull of the sets $G^+(P_k)$, $1\le
k\le p$. In this case, the set $E^\Phi=\{t^1,\ldots,t^N\}$ is
finite, it is just normals to the $(n-1)$-dimensional faces
$\Gamma_j(P)$ of the polyhedron situated outside the coordinate
planes, with  the condition $\Phi(t^j)=1$. The measure $\aoP$
charges $t^j$ with the volume of the convex hull $G_j^+(P)$ of
 the corresponding face $\Gamma_j(P)$ and $0$,
so
$$
\delta(dd^cu,\log|P|)\le\delta(dd^cu,\Psi_{\log|P|,0})=n!\,\sum_{1\le
j\le N}\psi_{u,0}(t^j)\,\Vol((G_j^+(P)).
$$

\vskip 0.5cm

Tek/Nat

H\o gskolen i Stavanger

POB 2557 Ullandhaug

N-4091 Stavanger

Norway

\vskip0.1cm

{\sc E-mail}: alexander.rashkovskii@tn.his.no

\end{document}